\documentclass[11pt]{article} 
\usepackage{amsmath}
\usepackage{amsfonts}
\usepackage{graphicx}
\usepackage{color}
\usepackage[latin1]{inputenc}
\usepackage[T1]{fontenc}
\usepackage[french]{babel}
\makeatletter
\def\@sect#1#2#3#4#5#6[#7]#8{%
  \ifnum #2>\c@secnumdepth
    \let\@svsec\@empty
  \else
    \refstepcounter{#1}%
    \protected@edef\@svsec{\@seccntformat{#1}\relax}%
  \fi
  \@tempskipa #5\relax
  \ifdim \@tempskipa>\z@
    \begingroup
      #6{%
        \@hangfrom{\hskip #3\relax\@svsec}%
          \interlinepenalty \@M #8\@@par}%
    \endgroup
    \csname #1mark\endcsname{#7}%
    \addcontentsline{toc}{#1}{%
      \ifnum #2>\c@secnumdepth \else
        \protect\numberline{\csname the#1\endcsname.}%
      \fi
      #7}%
  \else
    \def\@svsechd{%
      #6{\hskip #3\relax
      \@svsec #8}%
      \csname #1mark\endcsname{#7}%
      \addcontentsline{toc}{#1}{%
        \ifnum #2>\c@secnumdepth \else
          \protect\numberline{\csname the#1\endcsname.}%
        \fi
        #7}}%
  \fi
  \@xsect{#5}}
\def\@seccntformat#1{\csname the#1\endcsname.\quad}
\makeatother
\newtheorem{theo}[equation]{Th\'eor\`eme}
\newtheorem{fait}[equation]{Fait}
\newtheorem{lem}[equation]{Lemme}
\newtheorem{question}[equation]{Question}
\newtheorem{df}[equation]{D\'efinition}
\newtheorem{proposition}[equation]{Proposition}

\newenvironment{remarque}{
\refstepcounter{equation}\trivlist%
\item[\hskip \labelsep{\bfseries Remarque \theequation.\ }]}%
{\endtrivlist}%

\makeatletter
\renewcommand\theequation{\thesection.\arabic{equation}}
\@addtoreset{equation}{section}
\makeatother
\newcommand{\carrenoir}{\rule{0.5em}{0.5em}}
\makeatletter
\newenvironment{demo}[1][\@empty]{\textbf{D\'emonstration~%
\ifx\@empty#1:\else #1~:\fi~}}
{\hfill\carrenoir\nolinebreak\vspace{2mm}}
\makeatother
\newcommand{\oper}[2]{\newcommand{#1}{\mathop{\mathrm{#2}}\nolimits} }
\oper{\Vol}{Vol}
\newcommand{\R}{\mathbb R}
\newcommand{\N}{\mathbb N}
\oper{\SO}{SO}
\oper{\dimension}{dim}
\newcommand{\de}{\mathrm{ d }}

\title{Construction de valeurs propres doubles du laplacien de Hodge-de~Rham}
\author{Pierre Jammes}
\date{}
\begin{document}
\maketitle
{\small 
\textsc{Résumé.---}
Sur toute variété de dimension au moins 3, on construit une métrique
telle que la première valeur propre non nulle du laplacien agissant sur 
les $p$-formes différentielles soit double. On en déduit qu'on peut
prescrire le volume et le début du spectre du laplacien de Hodge-de~Rham
avec multiplicité 1 ou 2.

Mots-clefs : formes différentielles, laplacien de Hodge-de~Rham, 
multiplicité de valeurs propres.

\medskip
\textsc{Abstract.---}
On any compact manifold of dimension greater than 3, we exhib a metric
whose first positive eigenvalue for the Laplacian acting on $p$-form is of 
multiplicity 2. As a corollary, we prescribe the volume and any finite
part of the spectrum of the Hodge Laplacian with multiplicity 1 or 2.


}

\section{Introduction}
Y.~Colin de Verdière a montré dans \cite{cdv86} que pour toute variété
riemannienne compacte $M$ de dimension supérieure ou égale à~3 et tout
entier $N\geq1$, il existe une métrique sur $M$ telle que la multiplicité
de la première valeur propre du laplacien agissant sur les fonctions de $M$
soit égale à $N$, et a généralisé ce résultat en montrant qu'on peut
en fait prescrire toute partie finie du spectre du laplacien, la
multiplicité des valeurs propres pouvant être choisie arbitrairement (voir 
\cite{cdv87}).

 Le problème de la multiplicité des valeurs propres a aussi été étudié pour
des opérateurs de Schrödinger (\cite{cdvt93}, \cite{bcc98}), et les résultats 
de prescription de spectre a été adaptés par 
P.~Guérini au laplacien de Hodge-de~Rham ---~agissant sur les formes 
différentielles~--- dans \cite{gu04} et par M.~Dahl à l'opérateur de Dirac 
(\cite{da05}). Mais pour ces deux derniers opérateurs les valeurs propres 
prescrites sont simples (en se restreignant par exemple aux formes coexactes 
de degré fixé en ce qui concerne le laplacien de Hodge-de~Rham); le problème 
de créer des valeurs propres multiples sur une variété quelconque dans ces 
deux cas n'a pas encore été résolue.

 Le but de cet article est d'apporter un début de réponse à cette question
en expliquant comment construire des métriques telles qu'une valeur
propre du laplacien de Hodge-de~Rham soit de multiplicité~2. La
construction assure que cette multiplicité possède une
certaine stabilité, ce qui permet finalement de prescrire le début du
spectre avec multiplicité~1 ou 2.

 Comme dans \cite{ja06b}, on ne cherchera à prescrire
le spectre du laplacien qu'en restriction aux formes coexactes. En effet,
si on note 
\begin{equation}
0<\mu_{p,1}(M,g)\leq\mu_{p,2}(M,g)\leq\ldots
\end{equation} 
les valeurs propres du laplacien agissant sur les $p$-formes coexactes
de $M$, le spectre non nul du laplacien agissant sur l'ensemble des 
$p$-formes est $(\mu_{p-1,i}(M,g))_{i\geq1}\cup(\mu_{p,i}(M,g)_{i\geq1})$,
et la multiplicité de la valeur propre nulle, si elle est existe, est
le $p$-ième nombre de Betti de $M$. Comme la dualité de Hodge impose
que $\mu_{p,i}(M,g)=\mu_{n-p-1,i}(M,g)$, on peut se restreindre aux degrés
$p\leq[(n-1)/2]$.

\begin{theo}\label{intro:th1}
Soit $M$ une variété compacte, connexe et orientable de dimension $n=2k+1$
ou $2k+2$ où $k\in\N^*$, $V$ un réel strictement positif et $N\geq1$ un 
entier. On se donne pour tout entier $p\in\{1,\ldots,k\}$ une suite de réels 
$0<\nu_{p,1}\leq\nu_{p,2}\leq\ldots\leq\nu_{p,N}$, chaque valeur 
apparaissant au plus deux fois pour $p$ donné.

Il existe une métrique $g$ sur $M$ telle que
\begin{itemize}
\item $\mu_{p,i}(M,g)=\nu_{p,i}$ pour tout $i\leq N$ et $p\in\{1,\ldots,k\}$;
\item $\Vol(M,g)=V$.
\end{itemize}
\end{theo}

La construction de valeurs propres doubles s'appuiera sur l'apparition
d'un phénomène spectral lié à la présence de multiplicité, baptisé 
«~diabolo~» dans \cite{bw84}, dont on peut par exemple trouver la 
description dans \cite{ar76} (appendice~10) et \cite{cdv98} (chapitre~5).

Dans \cite{cdv88}, Y.~Colin de Verdière formalise une notion de 
transversalité remontant à Arnol'd et définit deux propriétés de 
transversalité ---~resp. forte et faible, que nous rappellerons
dans la section~\ref{stab}~--- pour les valeurs propres
multiples d'une forme quadratique. Il remarque dans \cite{cdv98}
que dans le cas où l'hypothèse forte est vérifiée pour une multiplicité~2, 
on peut montrer que cette multiplicité est stable mettant en évidence
un point diabolo. D'une certaine manière, notre démarche sera inverse~: 
on va construire un point diabolo sans faire intervenir d'hypothèse de 
transversalité, et on pourra vérifier \emph{a posteriori} que la valeur 
propre double vérifie l'hypothèse faible de transversalité.

Il serait bien sûr intéressant de savoir si on peut construire des
valeurs propres de multiplicité plus grande sur toute variété.
La technique que nous utiliserons est cependant
spécifique à la multiplicité~2 et ne semble pas pouvoir s'adapter
à des multiplicités plus grandes (voir remarque~\ref{stab:rq}).
Il convient toutefois de remarquer que sur certaines variétés, 
la multiplicité de la première valeur propre peut 
être arbitrairement grande:
\begin{theo}\label{intro:th2}
Pour tout entier $n\geq4$, et tout $1\leq p<n/2$, il existe une variété
$M$ de dimension $n$ telle que pour tout entier $k\geq1$, il existe une 
métrique $g$ sur $M$ telle que $\mu_{p,1}(M,g)$ soit de multiplicité au 
moins $k$.
\end{theo}
 En général, on ne peut donc pas majorer la multiplicité de la première 
propre en fonction de la topologie comme c'est le cas sur les surfaces.
Mais les exemples du théorème~\ref{intro:th2} ont une topologie particulière
(variétés produits) et on ne contrôle la multiplicité de la première valeur
propre que pour certains degrés qu'on ne peut pas choisir indépendamment 
de la topologie. En ce sens, le 
théorème~\ref{intro:th1} qui prescrit les premières valeurs propres avec 
multiplicité 1 ou 2 pour tous les degrés simultanément et sur n'importe 
quelle variété est beaucoup plus précis.  

 Ces résultats permettent de mieux cerner le problème de la multiplicité
des valeurs propres du laplacien de Hodge-de~Rham et de dégager quelques
questions qui restent en suspens, par exemple:
\begin{question}
La multiplicité de la première valeur propre de la sphère peut-elle être
arbitrairement grande quel que soit le degré ? 
\end{question}
\begin{question}
Si $M$ est une variété de dimension~3, la multiplicité de $\mu_{1,1}(M,g)$ 
est-elle nécessairement bornée ? Si oui, comment varie la multiplicité 
maximale de $\mu_{1,1}(M,g)$ en fonction de la topologie ?
\end{question} 

 La section~\ref{cv} sera consacrée au rappel des outils techniques
que nous utiliserons. Dans les sections~\ref{diabolo} et \ref{stab},
nous expliquerons comment construire une valeur propre double, et
pourquoi sa multiplicité est stable. Enfin, dans les sections~\ref{presc}
et \ref{exemple}, nous démontrerons les théorèmes~\ref{intro:th1}
et \ref{intro:th2}.

\section{Convergence de valeurs propres et d'espaces propres}\label{cv}
Nous allons rappeler ici les outils techniques qui vont intervenir dans
la construction de valeurs propres doubles. Le premier est le résultat
de convergence de valeurs propres et d'espaces propres obtenu 
par C.~Anné et B.~Colbois dans \cite{ac95} pour les variétés compactes
reliées par des anses fines : on se donne une famille finie de variétés
compactes $(M_j,g_j)_{i=1}^K$ qu'on relie entre elles par des anses fines,
isométriques au produit d'une sphère $(S^{n-1},\varepsilon^2 
g_{\mathrm{can}})$ par un intervalle. En notant 
$(\tilde M,g_\varepsilon)$ la variété obtenue, qui est difféomorphe
à $M_1\#M_2\#\ldots\#M_K$, on a alors :
\begin{theo}\label{global:th1}
Si, pour $p\in\{1,\cdots,n-1\}$, on note $\mu'_{p,1}\leq\mu'_{p,2}\leq\ldots$
la réunion des spectres $(\mu_{p,i}(M_j,g_j))_{i,j}$, on a pour tout 
$i\in\N^*$
$$\lim_{\varepsilon\to0}\mu_{p,i}(\tilde M,g_\varepsilon)=\mu'_{p,i},$$ 
et il y a convergence des espaces spectraux.
\end{theo}
Ce théorème a déjà été utilisé pour prescrire le spectre du laplacien
de Hodge-de~Rham sans multiplicité dans \cite{gu04} et \cite{ja06b}. 

Pour obtenir de la multiplicité, on aura besoin d'un contrôle sur la
vitesse de convergence des valeurs propres et des espaces propres. Ce
contrôle découle de la proposition 3.10 ainsi et des corollaires
3.11 et 3.12 de \cite{ac95} (voir aussi les propositions 1 et 2 de 
\cite{an90}. Nous renvoyons aussi à \cite{an90} pour une définition
précise de la convergence des espaces propres). Si $I$ est un intervalle
de $\R^+$, $E_I^\varepsilon$ l'espace engendré par les $p$-formes
propres de $(\tilde M,g_\varepsilon)$ de valeur propre contenue
dans $I$ et $E_I$ l'espace engendré par les $p$-formes propres 
des $(M_j,g_j)$ de valeur propre contenue dans $I$, la distance entre
$E_I$ et $E_I^\varepsilon$ vérifie $d(E_I,E_I^\varepsilon)\leq 
C\tau(\varepsilon)$ avec $\lim_{\varepsilon\to0}\tau(\varepsilon)=0$
et où $C$ est une constante dépendant des bornes de $I$ et de leur
distance au spectre $(\mu'_{p,i})_i$. En particulier, si deux 
valeurs propres de $(\tilde M,g_\varepsilon)$ sont proches l'une
de l'autre mais assez éloignées du reste du spectre, il est difficile de
localiser les formes propres mais la somme des deux espaces propres
converge rapidement.

Ces résultats seront appliqués à la variété obtenue en attachant à $M$
des sphères munies de métriques bien choisies. Plus précisément, 
on utilisera les «~haltères de Cheeger généralisées~» définies par 
P.~Guérini dans \cite{gu04} : pour tous entiers $n\geq3$ et $1\leq p\leq n-1$ 
et tout réel $u>0$ petit, 
on considère le domaine $\Omega_{p,u}$ de $\R^{n+1}$ formé par la réunion 
d'un $\frac15$-voisinage tubulaire de la sphère unité $S^p\subset 
R^{p+1}\times\{0\}\subset\R^{n+1}$ et du produit de boules 
$B^{p+1}(0,1)\times B^{n-p}(0,u)\subset\R^{p+1}\times\R^{n-p}$.
Après lissage de son bord, le domaine $\Omega_{p,u}$ est difféomorphe
à une boule $B^{n+1}$ et $\partial\Omega_{p,u}$ est une sphère (voir la 
section 2.1 de \cite{gu04} pour les détails de la construction). 
La famille de métriques induites sur le bord a la propriété de produire 
une petite valeur propre :
\begin{proposition}[\cite{gu04}]\label{cv:prop}
La famille $\partial\Omega_{p,u}$ vérifie
$$\lim_{u\to0}\mu_{p,1}(\partial\Omega_{p,u})=0$$
et il existe des constantes $c(n)>0$ et $v(n)>0$ telle que
$$\mu_{p,2}(\partial\Omega_{p,u})\geq c\textrm{ et }
\mu_{q,1}(\partial\Omega_{p,u})\geq c$$
pour tout $u$ et tout $q\neq p$, et
$$\Vol(\partial\Omega_{p,u})<v$$
pour tout $u$.
\end{proposition}

On aura besoin en outre de certaines propriétés de symétrie de
$\partial\Omega_{p,u}$. Par construction, le domaine $\Omega_{p,u}$
est invariant sous les actions de $\SO_{p+1}(\R)$ sur les $(p+1)$ premières
coordonnées de $\R^{n+1}$ et de $\SO_{n-p}(\R)$ sur les $(n-p)$
dernières coordonnées. On sait alors que si une valeur propre du laplacien
est petite par rapport au diamètre maximal des orbites, les formes
propres correspondantes sont invariantes (voir \cite{ja04}, théorème 1.10 
et lemme 4.1). Notons $\omega_{p,u}$ la forme propre de $\partial\Omega_{p,u}$
associée à la valeur propre $\mu_{1,p}(\partial\Omega_{p,u})$. Si
$u$ est suffisamment petit, la forme $\omega_{p,u}$ est donc invariante
par $\SO_{p+1}(\R)$ et $\SO_{n-p}(\R)$. En particulier, sa restriction
à la sphère $\partial\Omega_{p,u}\cap\R^{p+1}\times\{0\}$ est nécessairement,
à un facteur scalaire non nul près, la forme volume canonique.
Si on note $\Upsilon$ l'isométrie de $\R^{n+1}$ définie par
$\Upsilon(x_1,\ldots,x_{n+1})=(x_1,\ldots,x_p,-x_{p+1},-x_{p+2},x_{p+3},
\ldots,x_{n+1})$, 
on en déduit le 
\begin{fait}\label{cv:fait}
La sphère $\partial\Omega_{p,u}$ est invariante par $\Upsilon$ et 
$\Upsilon^*\omega_{p,u}=-\omega_{p,u}$.
\end{fait}
Cette propriété de la forme $\omega_{p,u}$ sera un élément clef de la 
construction de valeurs propres doubles.

\section{Création d'une valeur propre double}\label{diabolo}

On va maintenant utiliser les outils présentés dans la section précédente 
pour construire une valeur propre double.

\begin{figure}[h]
\begin{center}
\begin{picture}(0,0)%
\includegraphics{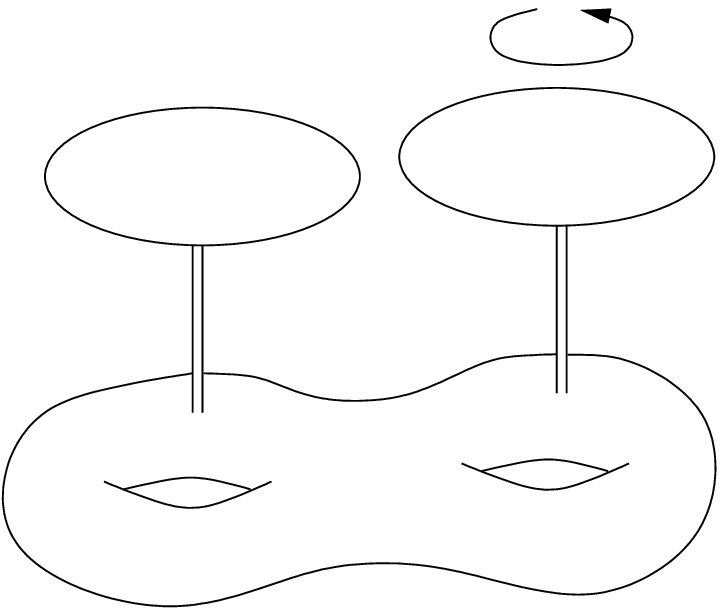}%
\end{picture}%
\setlength{\unitlength}{4144sp}%
\begingroup\makeatletter\ifx\SetFigFont\undefined%
\gdef\SetFigFont#1#2#3#4#5{%
  \reset@font\fontsize{#1}{#2pt}%
  \fontfamily{#3}\fontseries{#4}\fontshape{#5}%
  \selectfont}%
\fi\endgroup%
\begin{picture}(3283,2844)(1686,-2623)
\put(4636, 74){\makebox(0,0)[lb]{\smash{{\SetFigFont{12}{14.4}{\rmdefault}{\mddefault}{\updefault}{\color[rgb]{0,0,0}$\theta$}%
}}}}
\put(4006,-601){\makebox(0,0)[lb]{\smash{{\SetFigFont{12}{14.4}{\rmdefault}{\mddefault}{\updefault}{\color[rgb]{0,0,0}$\omega_2,\lambda_2$}%
}}}}
\put(2341,-646){\makebox(0,0)[lb]{\smash{{\SetFigFont{12}{14.4}{\rmdefault}{\mddefault}{\updefault}{\color[rgb]{0,0,0}$\omega_1,\lambda_1$}%
}}}}
\end{picture}%
\end{center}
\caption{\label{diabolo:anses}}
\end{figure}

 On se donne une variété riemannienne $(M,g)$ de dimension $n\geq3$ 
quelconque et un intervalle ouvert $I$ ne rencontrant pas le spectre
du laplacien agissant sur les $p$-formes coexactes de $(M,g)$, et on va 
créer une 
valeur propre double dans cet intervalle. En appliquant des homothéties
aux sphères $\partial\Omega_{p,u}$ de la proposition \ref{cv:prop},
on obtient sur la sphère $S^n$ deux métriques $g_1$ et $g_2$, ayant
une valeur propres $\lambda_i$, $i=1,2$ dans l'intervalle $I$ de
forme propre $\omega_i$, toutes les autres valeurs propres étant plus
grandes que $I$. En attachant les sphères $(S^n,g_1)$ et $(S^n,g_2)$
à $(M,g)$ par des anses fines de rayon $\varepsilon$, on obtient la variété 
de la figure \ref{diabolo:anses}.

 On prend soin de choisir comme point d'attache des anses sur les 
sphères un point fixe de l'isométrie $\Upsilon$, le recollement entre
l'anse et la sphère étant alors déterminé à une isométrie de $S^{n-1}$ près.
On fixe cette isométrie arbitrairement sur la 1\iere~sphère, et on
se laisse la liberté de faire varier cette isométrie dans un 
sous-groupe de $\SO(n)$ isomorphe à $S^1$ pour la 2\ieme~sphère, en notant
$\theta$ le paramètre naturel sur le cercle. La métrique ainsi définie
sur $M$ est paramétrée par $\lambda_1$, $\lambda_2$, $\varepsilon$ et
$\theta$. On peut en outre remarquer qu'on ne modifie pas cette
métrique en remplaçant $\theta$ par $\theta+\pi$, du fait de la
symétrie de la métrique $g_2$ (cf. fait \ref{cv:fait}).

\begin{figure}[h]
\begin{center}
\begin{picture}(0,0)%
\includegraphics{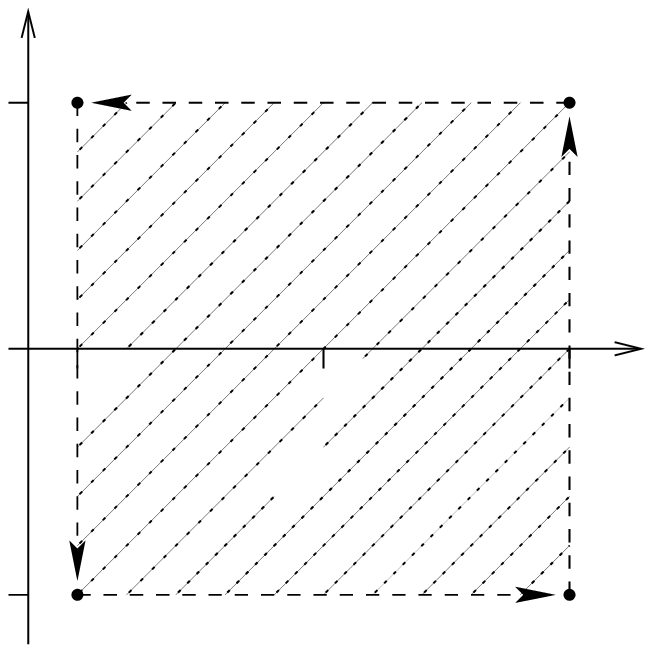}%
\end{picture}%
\setlength{\unitlength}{4144sp}%
\begingroup\makeatletter\ifx\SetFigFont\undefined%
\gdef\SetFigFont#1#2#3#4#5{%
  \reset@font\fontsize{#1}{#2pt}%
  \fontfamily{#3}\fontseries{#4}\fontshape{#5}%
  \selectfont}%
\fi\endgroup%
\begin{picture}(3540,3063)(481,-2605)
\put(856,299){\makebox(0,0)[lb]{\smash{{\SetFigFont{12}{14.4}{\rmdefault}{\mddefault}{\updefault}{\color[rgb]{0,0,0}$\lambda_2$}%
}}}}
\put(811,-1186){\makebox(0,0)[lb]{\smash{{\SetFigFont{12}{14.4}{\rmdefault}{\mddefault}{\updefault}{\color[rgb]{0,0,0}$\lambda_1$}%
}}}}
\put(496,-2311){\makebox(0,0)[lb]{\smash{{\SetFigFont{12}{14.4}{\rmdefault}{\mddefault}{\updefault}{\color[rgb]{0,0,0}$\lambda_1-\eta$}%
}}}}
\put(2476,-1366){\makebox(0,0)[lb]{\smash{{\SetFigFont{12}{14.4}{\rmdefault}{\mddefault}{\updefault}{\color[rgb]{0,0,0}$\frac\pi2$}%
}}}}
\put(1396,-1321){\makebox(0,0)[lb]{\smash{{\SetFigFont{12}{14.4}{\rmdefault}{\mddefault}{\updefault}{\color[rgb]{0,0,0}$0$}%
}}}}
\put(3646,-1321){\makebox(0,0)[lb]{\smash{{\SetFigFont{12}{14.4}{\rmdefault}{\mddefault}{\updefault}{\color[rgb]{0,0,0}$\pi$}%
}}}}
\put(4006,-1321){\makebox(0,0)[lb]{\smash{{\SetFigFont{12}{14.4}{\rmdefault}{\mddefault}{\updefault}{\color[rgb]{0,0,0}$\theta$}%
}}}}
\put(1306,-2536){\makebox(0,0)[lb]{\smash{{\SetFigFont{12}{14.4}{\rmdefault}{\mddefault}{\updefault}{\color[rgb]{0,0,0}$t_1$}%
}}}}
\put(496,-61){\makebox(0,0)[lb]{\smash{{\SetFigFont{12}{14.4}{\rmdefault}{\mddefault}{\updefault}{\color[rgb]{0,0,0}$\lambda_1+\eta$}%
}}}}
\put(1306, 74){\makebox(0,0)[lb]{\smash{{\SetFigFont{12}{14.4}{\rmdefault}{\mddefault}{\updefault}{\color[rgb]{0,0,0}$t_0$}%
}}}}
\put(3556, 74){\makebox(0,0)[lb]{\smash{{\SetFigFont{12}{14.4}{\rmdefault}{\mddefault}{\updefault}{\color[rgb]{0,0,0}$t_3$}%
}}}}
\put(3556,-2536){\makebox(0,0)[lb]{\smash{{\SetFigFont{12}{14.4}{\rmdefault}{\mddefault}{\updefault}{\color[rgb]{0,0,0}$t_2$}%
}}}}
\put(2296,-1816){\makebox(0,0)[lb]{\smash{{\SetFigFont{14}{16.8}{\rmdefault}{\mddefault}{\updefault}{\color[rgb]{0,0,0}$D$}%
}}}}
\end{picture}%
\end{center}
\caption{\label{diabolo:param}}
\end{figure}

On fixe $\lambda_1$ dans l'intervalle $I$ et on choisit une réel $\eta>0$
tel que $\lambda_1\pm\eta\in I$. On va créer de la multiplicité en 
faisant varier $\lambda_2$ dans l'intervalle $[\lambda_1-\eta,
\lambda_1+\eta]$ et $\theta$ dans l'intervalle $[0,\pi]$. On fixe
$\varepsilon$ suffisamment petit pour que $M$ n'ait que deux valeurs
propres $\mu_1$ et $\mu_2$ dans l'intervalle $I$, proches de $\lambda_1$ et
$\lambda_2$, et que la somme de leurs espaces propres soit proche de 
l'espace engendré par $\omega_1$ et $\omega_2$. Quand les valeurs propres 
$\mu_1$ et $\mu_2$ sont distinctes, on posera $\mu_1<\mu_2$ et on notera
$\varphi_1$ et $\varphi_2$ leurs formes propres respectives. On peut
en outre choisir $\varepsilon$ de sorte que si $\lambda_2=\lambda_1+\eta$
(resp. $\lambda_1-\eta$), $\varphi_1$ est proche de $\omega_1$ (resp.
$\omega_2$). Au final, on ne fera varier la métrique que dans le
domaine $D$ de dimension~2 représenté sur la figure \ref{diabolo:param},
et c'est dans ce domaine qu'on va trouver une valeur propre double.

\begin{proposition}\label{diabolo:prop}
Il existe un point de $D$ pour lequel $\mu_1=\mu_2$.
\end{proposition}
\begin{demo}
L'idée est de considérer la famille à un paramètre de métriques $(g_t)$
obtenue en suivant le bord du domaine $D$ et de suivre l'évolution des
$\mu_i$ et $\varphi_i$. On note $t_i$ les valeurs prises par $t$ à 
chaque coin de $D$ (voir figure~\ref{diabolo:param}). Pour toute les
métriques paramétrées par $D$, l'espace engendré par les $\varphi_i$
est proche de celui engendré par les $\omega_i$, il existe donc
un isométrie naturelle entre les deux qui permet de les identifier et
de représenter les $\varphi_i$ et $\omega_i$ dans un même plan, comme
dans la figure~\ref{diabolo:vp}.
\begin{figure}[h]
\begin{center}
\begin{minipage}[b]{.46\linewidth}
\begin{picture}(0,0)%
\includegraphics{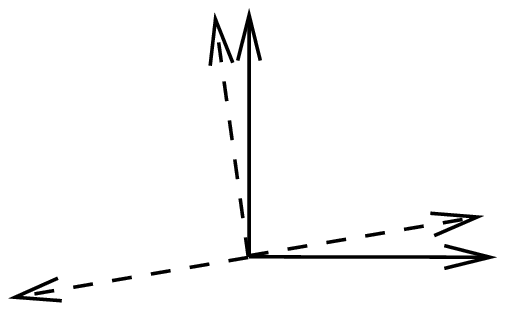}%
\end{picture}%
\setlength{\unitlength}{4144sp}%
\begingroup\makeatletter\ifx\SetFigFont\undefined%
\gdef\SetFigFont#1#2#3#4#5{%
  \reset@font\fontsize{#1}{#2pt}%
  \fontfamily{#3}\fontseries{#4}\fontshape{#5}%
  \selectfont}%
\fi\endgroup%
\begin{picture}(2346,1694)(911,-1846)
\put(1018,-1781){\makebox(0,0)[lb]{\smash{{\SetFigFont{10}{12.0}{\rmdefault}{\mddefault}{\updefault}{\color[rgb]{0,0,0}$\varphi_1(t_3)$}%
}}}}
\put(3040,-1586){\makebox(0,0)[lb]{\smash{{\SetFigFont{10}{12.0}{\rmdefault}{\mddefault}{\updefault}{\color[rgb]{0,0,0}$\omega_1$}%
}}}}
\put(3040,-1120){\makebox(0,0)[lb]{\smash{{\SetFigFont{10}{12.0}{\rmdefault}{\mddefault}{\updefault}{\color[rgb]{0,0,0}$\varphi_1(t_0)$}%
}}}}
\put(1446,-342){\makebox(0,0)[lb]{\smash{{\SetFigFont{10}{12.0}{\rmdefault}{\mddefault}{\updefault}{\color[rgb]{0,0,0}$\varphi_1(t_1)$}%
}}}}
\put(2146,-303){\makebox(0,0)[lb]{\smash{{\SetFigFont{10}{12.0}{\rmdefault}{\mddefault}{\updefault}{\color[rgb]{0,0,0}$\omega_2$}%
}}}}
\end{picture}%
\caption{\label{diabolo:vp}} 
\end{minipage}\hfill
\begin{minipage}[b]{.46\linewidth}
\begin{picture}(0,0)%
\includegraphics{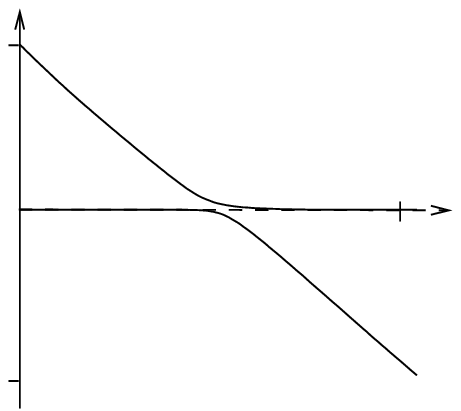}%
\end{picture}%
\setlength{\unitlength}{4144sp}%
\begingroup\makeatletter\ifx\SetFigFont\undefined%
\gdef\SetFigFont#1#2#3#4#5{%
  \reset@font\fontsize{#1}{#2pt}%
  \fontfamily{#3}\fontseries{#4}\fontshape{#5}%
  \selectfont}%
\fi\endgroup%
\begin{picture}(2622,1869)(211,-1423)
\put(2476,-1051){\makebox(0,0)[lb]{\smash{{\SetFigFont{12}{14.4}{\rmdefault}{\mddefault}{\updefault}{\color[rgb]{0,0,0}$\mu_1$}%
}}}}
\put(1081, 74){\makebox(0,0)[lb]{\smash{{\SetFigFont{12}{14.4}{\rmdefault}{\mddefault}{\updefault}{\color[rgb]{0,0,0}$\mu_2$}%
}}}}
\put(856,-646){\makebox(0,0)[lb]{\smash{{\SetFigFont{12}{14.4}{\rmdefault}{\mddefault}{\updefault}{\color[rgb]{0,0,0}$t_0$}%
}}}}
\put(2521,-691){\makebox(0,0)[lb]{\smash{{\SetFigFont{12}{14.4}{\rmdefault}{\mddefault}{\updefault}{\color[rgb]{0,0,0}$t_1$}%
}}}}
\put(536,-581){\makebox(0,0)[lb]{\smash{{\SetFigFont{12}{14.4}{\rmdefault}{\mddefault}{\updefault}{\color[rgb]{0,0,0}$\lambda_1$}%
}}}}
\put(226,219){\makebox(0,0)[lb]{\smash{{\SetFigFont{12}{14.4}{\rmdefault}{\mddefault}{\updefault}{\color[rgb]{0,0,0}$\lambda_1+\eta$}%
}}}}
\put(226,-1316){\makebox(0,0)[lb]{\smash{{\SetFigFont{12}{14.4}{\rmdefault}{\mddefault}{\updefault}{\color[rgb]{0,0,0}$\lambda_1-\eta$}%
}}}}
\put(2706,-436){\makebox(0,0)[lb]{\smash{{\SetFigFont{12}{14.4}{\rmdefault}{\mddefault}{\updefault}{\color[rgb]{0,0,0}$t$}%
}}}}
\end{picture}%
\caption{\label{diabolo:graphe}} 
\end{minipage}
\end{center}
\end{figure}

Quand $t=t_0$ on a $\lambda_2=\lambda_1+\eta$, donc $\mu_1$ et $\varphi_1$
sont proches de $\lambda_1$ et $\omega_1$ respectivement (voir figures~%
\ref{diabolo:vp} et \ref{diabolo:graphe}).

 Quand on passe de $t_0$ à $t_1$, il est possible ---~mais improbable~---
qu'on croise un point pour lequel $\mu_1=\mu_2$. Si c'est le cas, la
proposition est démontrée. Sinon, les deux valeurs propres restent 
distinctes et évoluent
comme sur la figure~\ref{diabolo:graphe}~: pour $t=t_1$, les valeurs 
propres $\mu_1$ et $\mu_2$ sont proches de $\lambda_2=\lambda_1-\eta$ et 
$\lambda_1$ respectivement. La forme $\varphi_1$ s'est donc déplacée et
est maintenant
proche de $\omega_2$ comme sur la figure~\ref{diabolo:vp} (quitte à changer 
la convention de signe sur $\omega_2$). 

 On passe ensuite de $t_1$ à $t_2$, c'est-à-dire qu'on fait tourner
la 2\ieme{} sphère d'un demi-tour autour du point d'attache à l'anse
sans modifier $\lambda_1$ et $\lambda_2$. La forme $\varphi_1$ reste
donc proche de $\omega_2$. Il faut noter que les métriques $g_{t_1}$
et $g_{t_2}$ sont isométriques, et donc qu'on peut identifier les espaces 
propres en $t_1$ et $t_2$, mais que $\omega_2$ et $\varphi_1$ ont changé 
de signe par rapport à cette identification (cf. fait \ref{cv:fait}).

 On va ensuite de $t_2$ à $t_3$ en faisant passer $\lambda_2$ de
$\lambda_1-\eta$ à $\lambda_1+\eta$. Par rapport au chemin $[t_1,t_0]$,
la différence est que $\theta$ a varié de $\pi$, c'est-à-dire qu'on
passe par des métriques isométriques à celle de $[t_1,t_0]$. Les
valeurs propres varient donc exactement comme dans la 
figure~\ref{diabolo:graphe}
en inversant l'axe temporel (en particulier, les valeurs propres
$\mu_1$ et $\mu_2$ restent distinctes). Les espaces propres subissent 
aussi l'évolution inverse de celle du parcours $t_0\to t_1$. Mais 
$\varphi_1$ a changé de signe entre $t_1$ et $t_2$, en $t_3$ 
c'est donc $-\varphi_1$ qui est proche de $\omega_1$. Les métriques
$g_{t_0}$ et $g_{t_3}$sont isométriques, mais la forme $\varphi_1$ a
changé de signe entre $t_0$ et $t_3$.

Entre $t_3$ et $t_0$, la forme $\varphi_1$ reste proche de $-\omega_1$
puisque $\lambda_1$ et $\lambda_2$ ne varient pas. Quand la métrique suit
le bord de $D$, la forme $\varphi_1$ associée à la valeur propre $\mu_1$
varie donc continûment jusqu'à prendre finalement, la valeur
opposée de sa valeur initiale.

 Considérons maintenant le fibré vectoriel trivial au dessus de $D$ de 
fibre $\R^2$, en identifiant la fibre à l'espace engendré par 
$\omega_1$ et $\omega_2$, c'est-à-dire au plan représenté par la 
figure~\ref{diabolo:vp}. Comme il n'y a pas de multiplicité en restriction
à $\partial D$, le sous-fibré en droite, de fibre $\R\varphi_1$ est
bien défini. On vient de montrer que ce fibré est non orientable.

Supposons maintenant que $\mu_1\neq\mu_2$ pour toutes les valeurs de
$\lambda_2$ et $\theta$ dans $D$. Le fibré en droite induit par
$\varphi_1$ est alors défini sur le domaine $D$ entier. Or, ce domaine
est contractile, donc ce fibré est trivial. En particulier, sa restriction
à $\partial D$ ne peut pas être non orientable. On en déduit par
l'absurde qu'il existe un point à l'intérieur de $D$ pour lequel
$\mu_1=\mu_2$.
\end{demo}

\section{Stabilité de la multiplicité}\label{stab}
 Commençons par rappeler la définition des l'hypothèses de transversalité
donnée par Y.~Colin de Verdière dans \cite{cdv88} pour les valeurs propres
multiples (voir aussi la section~3 de \cite{an90}).

 On se donne un intervalle $I$ de $\R$, un famille de métrique $g_a$,
où le paramètre $a$ prend ses valeurs dans un compact $K$ d'une variété,
que l'on supposera difféomorphe à une boule, 
et on suppose que pour une valeur $a_0$ du paramètre, on a une unique 
valeur propre $\lambda_0$ du 
laplacien dans $I$ pour la métrique $g_0$, de multiplicité
$n_0$, et d'espace propre $E_0$. Pour $a$ proche de $a_0$, la somme
$E_a$ des espaces propres de valeurs propres contenu dans $I$ est proche
de $E_0$, ce qui permet de les identifier par une isométrie naturelle (voir
\cite{cdv88} ou \cite{an90} pour sa définition explicite) et on note
$q_a\in\mathcal Q(E_0)$ la forme quadratique ainsi induite sur $E_0$ par
le laplacien. On note $\Phi:K\to\mathcal Q(E_0)$ l'application $a\mapsto q_a$
ainsi définie.
\begin{df}
La valeur propre $\lambda_0$ vérifie l'hypothèse forte (resp. faible) 
de transversalité si $\Phi$ est une submersion en $a=a_0$ (resp. $\Phi$ 
est essentielle en $a_0$).
\end{df}
Comme dans \cite{an90} on dira dans la suite de ce texte qu'une telle
valeur propre est \emph{fortement} (resp. \emph{faiblement}) \emph{stable}. 
La définition d'une application essentielle en $a_0$ utilisée dans 
\cite{cdv88} est que si pour tout $\Psi:K\to\mathcal Q(E_0)$ tel que
$\|\Psi-\Phi\|_\infty<\varepsilon$ on a $\Phi(a_0)\in\Psi(K)$ (nous
dirons ici qu'une telle application est \emph{métriquement essentielle}).
unous utiliserons une notion plus forte, que nous nommerons 
\emph{topologiquement essentielle}, qui est que la restriction de 
$\Phi$ à $\partial K$ n'est pas homotopiquement triviale dans $\mathcal Q(E_0)
\backslash\{\Phi(a_0)\}$. On peut vérifier que si $\Phi$ est topologiquement
essentielle en $a_0$ elle est aussi métriquement essentielle, et que
si deux applications sont topologiquement essentielles en $a_0$ et $b_0$ 
respectivement, alors leur produit est essentiel en $(a_0,b_0)$.

\begin{proposition}\label{stab:prop}
La valeur propre double de la proposition~\ref{diabolo:prop} est faiblement
stable.
\end{proposition}
\begin{demo}
Dans la section précédente, Le compact $K$ est le domaine $D$, l'espace
$E_0$ est l'espace engendré par $\omega_1$ et $\omega_2$ représenté par
la figure~\ref{diabolo:vp}. Pour normaliser le problème, on se ramène par 
homothétie à des formes quadratiques sur $E_0$ de trace fixée, c'est-à-dire 
qu'on pose $\mu_1+\mu_2=c^{te}$. Cet espace est de dimension~2,
on peut l'assimiler à un plan dont un seul point représente une forme
quadratique telle que $\mu_1=\mu_2$. Le reste du plan peut être
paramétré en coordonnées polaire par la différence entre les deux valeurs
propres et la direction de la droite propre correspondant par exemple à 
$\mu_2$ (attention: en faisant tourner la droite propre d'un demi-tour,
on revient à la forme quadratique initiale, c'est-à-dire qu'on a fait
une tour complet dans le plan des formes quadratiques). 
En faisant tourner les droites propres comme l'indique la 
figure~\ref{diabolo:vp}, on a en fait montré que $\Phi$ envoie $\partial D$
sur une courbe du plan qui entoure la forme quadratique $q_0$ de 
multiplicité~2, c'est-à-dire que l'indice de la courbe $\Phi(\partial D)$ 
par rapport à $q_0$ est non nul. L'application $\Phi$ est donc
topologiquement essentielle en un point de $D$ d'image $q_0$, et
la multiplicité est faiblement stable (on retrouve au passage 
l'existence de la valeur propre double). Dans le cas où la multiplicité
apparaît sur l'axe $\theta=0$, on al fait disparaître en changeant l'origine
due paramètre $\theta$ ou en perturbant la métrique $g$ sur la variété $M$. 
\end{demo}
\begin{remarque}
La stabilité forte d'une valeur propre multiple est une notion
différentielle, alors que la stabilité faible est seulement topologique.
Tous les arguments que nous avons utilisé dans cette section et la précédente
sont de nature purement topologique, il ne semble donc pas envisageable de 
montrer que la stabilité est forte par ces moyens.
\end{remarque}
\begin{remarque} La stabilité de la multiplicité donnée par la 
proposition~\ref{diabolo:prop} peut s'exprimer d'une autre manière: si on 
déforme légèrement le domaine $D$ dans l'espace des métriques de sorte que 
le fibré en droite au dessus du bord $\partial D$ reste non trivial, 
l'argument topologique utilisé reste valable, il y a donc toujours à 
l'intérieur du domaine un point donnant une valeur propre double.
\end{remarque}
\begin{remarque}\label{stab:rq}
Une condition nécessaire pour qu'une multiplicité soit stable est
que $\dimension K\geq\dim\mathcal Q(E_0)$ (ou $\dim\mathcal Q(E_0)$-1
si on se restreint à des formes quadratiques de trace fixée). En particulier,
la dimension de $K$ doit croître au moins quadratiquement par rapport à
la multiplicité. La technique que nous avons utilisée pour construire
une valeur propre double ne semble donc pas pouvoir se généraliser
à des multiplicités supérieures. En effet, si on attache un plus grand
nombre de sphères à la variété, le nombre de paramètres 
disponibles restera une fonction affine du nombre de sphères, même
en déplaçant le point d'attache sur la variété et en utilisant tous
les degrés de liberté à la jonction des anses.
\end{remarque}

On peut déduire des résultats obtenus jusqu'ici l'énoncé suivant, qui sera
la base de la démonstration du théorème~\ref{intro:th1}: 
\begin{lem}\label{stab:lem}
Soit $M$ une variété de dimension $n\geq3$, $p$ un entier tel que
$1\leq p\leq\frac{n-1}2$ et $C>\nu>0$ et $V>0$ trois réels. Il existe une
métrique $g$ sur $M$ telle que 
\begin{itemize}
\item la valeur propre $\mu_{p,1}(M,g)=\nu=\mu_{p,2}(M,g)$ est de 
multiplicité~2 et faiblement stable;
\item $\mu_{p,3}(M,g)>C$, et $\mu_{q,1}(M,g)>C$ pour $1\leq q\leq\frac{n-1}2$
et $q\neq p$;
\item $\Vol(M,g)<V$.
\end{itemize}
\end{lem}
\begin{demo} On commence par choisir sur $M$ une métrique $g$ telle que
$\Vol(M,g)<V/10$ et $\mu_{q,1}(M,g)>2C$ pour tout $q$. On applique ensuite
la construction de la section~\ref{diabolo} en considérant l'intervalle 
$I=]\frac\nu2,\frac{3\nu}2[$, en attachant des sphères
de volume inférieur à $V/10$ et telle que leurs valeurs propres qui ne sont 
pas dans $I$ soient plus grandes que $2C$ (la proposition~\ref{cv:prop}
le permet) et en choisissant le rayon $\varepsilon$ 
des anses de sorte que le volume total soit inférieur à $4V/10$. Les
propositions~\ref{diabolo:prop} et \ref{stab:prop} nous disent qu'il y a 
alors une valeur propre
$\mu_{p,1}(M,g)$ de multiplicité~2 faiblement stable dans l'intervalle $I$.
Une homothétie permet de ramener cette valeur propre à $\mu_{p,1}(M,g)=\nu$
et les conclusions du lemme sont vérifiées.
\end{demo}
\begin{remarque}\label{stab:simple}
Le même énoncé pour une valeur propre simple découle immédiatement de la 
proposition~\ref{cv:prop}, la stabilité étant trivialement vérifiée dans
ce cas.
\end{remarque}

\section{Prescription du spectre}\label{presc}
En s'appuyant sur le lemme~\ref{stab:lem} et en utilisant les techniques
habituelles de stabilité spectrale, on peut prescrire le début du 
spectre avec multiplicité~1 ou~2.

\begin{demo}[du théorème~\ref{intro:th1}]
On fixe un réel $\delta>0$ tel que
\begin{equation}
\delta<\inf_{p<k,\nu_{p,i}\neq\nu_{p,j}}\left\{\frac{|\nu_{p,i}-\nu_{p,j}|}2
\right\} \textrm{ et } \delta<\frac V2
\end{equation}
et une constante $C>\sup_{p,i}\nu_{p,i}$. On notera $v$ un paramètre
variant dans l'intervalle $[V-\delta,V+\delta]$.

Soit $m$ le nombre total de valeur propre à prescrire, tous degrés confondus
mais sans compter la multiplicité. Pour chaque $i\leq m$, on note $\nu$ la 
valeur propre correspondante, on se donne une sphère munie de la métrique 
donnée par le lemme~\ref{stab:lem} si la valeur propre est double et la 
remarque~\ref{stab:simple} si elle est simple, avec un volume $v_i$ inférieur
à $V/2m$: si on note $E_i$ l'espace propre correspondant et $q_i$ la forme
quadratique sur $E_i$ de valeur propre $\nu$, il existe un compact $K_i$,
une famille de métriques $g_{a_i}$ paramétrée par $a_i\in K_i$,
un point $\bar a_i\in K_i$ et une application 
$\Phi_i:K_i\to\mathcal Q(E_i)$ tels que $\Phi$ soit essentielle en $\bar a_i$
avec $\Phi_i(\bar a_i)=q_i$, $\Phi_i(a_i)$ étant la forme quadratique
induite sur $E_i$ par le laplacien pour une métrique $g_{a_i}$ sur la
sphère. On choisit $K_i$ suffisamment petit pour que les valeurs
propres de $\Phi_i(a_i)$ restent dans l'intervalle $I_i=
[\nu-\delta,\nu+\delta]$. Par définition de $\delta$ les intervalles
$I_i$ sont tous disjoints.

 On munit la variété $M$ d'une métrique $g$ telle que $\Vol(M,g)=
v-\sum_{i=1}^m v_i$ et $\mu_{p,1}(M,g)>C$ pour $1\leq p\leq n/2$ (c'est
possible selon \cite{gp95}) et on attache les $m$ sphères précédemment 
obtenues, munies d'une métrique $g_{a_i}$ avec $a_i\in K_i$, à $(M,g)$ par 
des anses de rayon $\varepsilon$. On obtient une variété difféomorphe à 
$M$, munie d'une métrique qu'on notera $g_\varepsilon$ et qui dépend
non seulement de $\varepsilon$ mais aussi des $a_i$.

La convergence du spectre et des espaces propres
donnée par le théorème~\ref{global:th1} permet de définir ainsi une
application
\begin{equation}
\Phi_\varepsilon:\left(\prod_{i=1}^kK_i\right)\times[V-\delta,V+\delta]
\to\left(\prod_{i=1}^k\mathcal Q(E_i)\right)\times\R,
\end{equation}
dont les premières composantes sont données par le spectre de 
$(M,g_\varepsilon)$ et la dernière étant le volume de $M$, qui converge
simplement ---~et donc uniformément d'après le théorème de Dini~---
vers l'application $\left(\prod_{i=1}^k\Phi_i\right)\times Id$ qui est 
topologiquement, donc métriquement, essentielle en 
$(\bar a_1,\ldots,\bar a_m,V)$. On en déduit que pour $\varepsilon$ 
suffisamment petit, le début du spectre et le volume souhaités sont donnés 
par un élément de l'image de $\Phi_\varepsilon$ et que les autres valeurs 
propres sont plus grandes que $C$.
\end{demo}
\section{Grandes multiplicités}\label{exemple}

 Pour finir, on va démontrer le théorème~\ref{intro:th2}. Étant donnés
$n\geq4$ et $1\leq p<n/2$, on considère une variété $N_1$ de dimension~3
et une variété $N_2$ de dimension~$n-3$ telle que $b_p(N_2)>0$, on pose 
$M=N_1\times N_2$ et on fixe un entier $k$ quelconque.  Comme la variété 
$N_1$ est de dimension~3, on peut la munir d'une métrique $g_1$ telle que 
la première valeur propre non nulle $\mu_{0,1}(N_1,g_1)$ du laplacien 
agissant sur les fonctions soit de multiplicité $k$ (cf.~\cite{cdv86}). 
Une particularité de la dimension~3 est qu'on peut choisir 
cette métrique de sorte que $\mu_{1,1}(N_1,g_1)>\mu_{0,1}(N_1,g_1)$.
En effet, les métriques qui interviennent dans la construction de
la valeur propre multiple sont obtenues par des déformations conformes
à partir d'une famille $(g_a)$ de métriques bien choisies, le volume
de $(N_1,g_1)$ pouvant être arbitrairement petit (voir la partie~4 de 
\cite{cdv86}). Or, selon \cite{ja06a},
si pour toute métrique $g$ sur $N_1$ on pose $C(g)=\inf_{\theta\in\Omega^1(N_1)}
\sup_{\de\zeta=0}(\|\de\theta\|_{\frac32}^2/\|\theta-\zeta\|_3^2)$
où $\|\cdot\|_p$ désigne la norme $L^p$, on a 
$\mu_{1,1}(N_1,\tilde g)\Vol(N_1,\tilde g)\geq C(g)$ pour toute métrique 
$\tilde g$ appartenant à la classe conforme de $g$. On peut en outre
choisir la famille $(g_a)$ telle qu'elle soit petite pour la distance
de Lipschitz. La constante $C(g)$ sera alors clairement minorée sur cette
famille, on a donc une minoration uniforme $\mu_{1,1}(N_1,\tilde g)
\Vol(N_1,\tilde g)\geq C>0$ pour toute métrique $\tilde g$ conforme à une
métrique $g_a$, et en particulier on peut choisir $g_1$ telle que 
$\mu_{1,1}(N_1,g_1)\geq C/\Vol(N_1,g_1)>\mu_{0,1}(N_1,g_1)$.

On munit ensuite $N_2$ d'une métrique 
$g_2$ telle que le spectre non nul de $(N_2,g_2)$ soit plus grand que 
$\mu_{0,1}(N_1,g_1)$. Si $M$ est muni de la métrique produit 
$g_M=g_1\oplus g_2$, la formule de Künneth nous dit que 
$\mu_{p,1}(M,g_M)=\mu_{0,1}(N_1,g_1)$, la multiplicité de $\mu_{p,1}(M,g_M)$ 
étant au moins égale à celle de $\mu_{0,1}(N_1,g_1)$, les produits des relevées 
des formes propres de $\mu_{0,1}(N_1,g_1)$ avec une $p$-forme harmonique 
de $(N_2,g_2)$ étant des formes propres de $\mu_{p,1}(M,g_M)$.

\begin{remarque}
La multiplicité $k$ étant stable sur $N_1$, on peut la transplanter sur
la somme connexe de $M$ avec une autre variété à l'aide du 
théorème~\ref{global:th1}.
\end{remarque}

\noindent Pierre \textsc{Jammes}\\
Université d'Avignon\\
laboratoire de mathématiques\\
33 rue Louis Pasteur\\
F-84000 Avignon\\
\texttt{Pierre.Jammes@univ-avignon.fr}
\end{document}